\documentclass[a4paper]{amsart}
\usepackage{amsmath,amssymb,amsfonts}
\usepackage{array}
\newcommand{\Z}{\mathbb{Z}}

\newtheorem{theorem}{Theorem}[section]
\newtheorem{corollary}[theorem]{Corollary}
\newenvironment{prf}[1]
{\begin{trivlist}%
\item {\bf Proof of \ref{#1}.}}
{\hspace*{\fill}$\Box$%
\end{trivlist}}
\title{Riemann surfaces with a large abelian group of automorphisms}
\thanks{This paper is part of the author's thesis. The
author is indebted with its advisor Rita Pardini for suggesting the
problem and for constant encouragment and advice.
}
\begin{document}
\maketitle\vspace{-.6cm}
\begin{center}
\mdseries\footnotesize
CLELIA LOMUTO\\[1mm]
\scshape\scriptsize
Dipartimento di Matematica "Ulisse Dini"
Viale Morgagni 67a\\
50134 FIRENZE, Italy.\\[.7cm]
\parbox{.79\textwidth}{%
\footnotesize\upshape 
{\scshape Abstract.} In this paper we classify all Riemann surfaces having a
  large abelian group of automorphisms, that is having an abelian group of
  automorphism of order strictly bigger then $4(g-1)$, where $g$
  denotes as usual the genus of the Riemann surface.}
\end{center}\vspace{3mm}
\section{Introduction}
\noindent In $1890$ Schwartz \cite{schwartz} proved that a compact Riemann
surface of genus $g\geq 2$ has a finite number of automorphisms.\\
More precisely, a compact Riemann surface  of genus $g\geq 2$ has at
most $84(g-1)$ automorphisms: this bound was given by Hurwitz as an
application of the now famous Riemann-Hurwitz formula \cite{hur} in $1893$.\\
This bound is sharp, as shown by the Klein quartic $\{x^3y+ y^3z+
z^3x=0 \}\subset \mathbb{P}^2$ studied by Klein \cite{klein} in 1879,
a genus $3$ curve with 168 automorphisms.\\
However, better bounds can be found with more restrictive assumptions.
In $1879$ Wiman \cite{wiman}
showed that a cyclic group of automorphisms of a Riemann 
surface cannot have order larger then $4g+2$.
In \cite{b4} Nakajima proved the upper bound $4g+4$ for an abelian group
of automorphisms.\\ 
In this paper we are interested in the problem of finding Riemann
surfaces with ``many'' automorphisms. More precisely,
given a compact Riemann surface $X$ of genus $g\geq 2$ and $G$ group of automorfisms of $X$, $G$ is called {\it large} if
$|G|>4(g-1)$.\\
In section $2$ we recall some known facts about group actions on
Riemann surfaces. In particular we recall  (cf.~\cite{b5}) that if 
$G$ is a large group of automorphisms of $X$ and $Y=X/G$, then
$Y$ has genus $g=0$ and using Riemann-Hurwitz formula one can give some
restriction on the number of critical values and on the ramification
indices of the projection map $X \to Y$. \\
In section $3$ we refine the above mentioned restrictions under the
assumption that the considered large group of automorphisms is
abelian: in theorem \ref{listaclelia} we reduce to a
list of $3$ infinite families and $11$ exceptional 
cases. The proof is based on some results of Nakajima \cite{b4} about
the stabilizer subgroups of the action of an abelian group of automorphisms.
\\
In section $4$, we construct all compact Riemann surfaces
with an abelian large group of automorphisms, showing that all cases
in the list in the previous section do occur. In each
case we give equations of a singular model of the curve. Finally
we prove also that such Riemann surfaces are hyperelliptic for genus
$g\geq 8$.\\
The basic tool used in the proof is the structure theorem for Galois
covering, due to Pardini \cite{b3}.
\begin{quote}
{\bf Notation:} From now on we will write 
Riemann surface $X$ for compact Riemann surface $X$ of genus $g(X)\geq2$.
\end{quote}
\section{Riemann surfaces having a large group of automorphisms}
\noindent In this section we briefly recall some facts about group actions on Riemann
surfaces. For further details, see \cite{b2}.\\
Let $G$ be a finite group acting holomorphically and effectively on a 
Riemann surface $X$, and let $P\in X$ be a fixed point. Then:\\
(a) the stabilizer subgroup $H_P$ is cyclic;\\
(b) there are only finitely many points in $X$ having nontrivial
stabilizer;\\
(c) there is a (unique) structure of Riemann surface on the quotient $X/G$ such
that the projection $\pi: X\to X/G$ is holomorphic; $\pi$ has degree $|G|$, and
for any point $P\in X$, $\textup{mult}_P(\pi)=|H_P|$.\\
Moreover, assume that $H_P$ is nontrivial and let $g$ be a
generator of the stabilizer subgroup $H_P$. Then there is a local
coordinate $z$ on $X$ centered at $P$ 
such that $g(z)=\lambda z$, where $\lambda$ is a primitive
$|H_P|^{th}$ root of unity.\\[3mm]
It follows from above that the quotient map $\pi:X\to X/G$ is a $|G|$-sheeted
branched cover of $X/G$: such a map is called Galois cover (with
Galois group $G$). \\
Let $Q_1,\ldots,Q_s$ be the critical values
of $\pi$, and let us associate to each $Q_i$ its {\it ramification
  index} $v_i=|H_P|$, where $H_P$ is the stabilizer subgroup of any
$P\in\pi^{-1}(Q_i)$. Note that the ramification index of a critical
value is well-defined since if $P$ and $P'$ are points of $X$ with
$\pi(P)=\pi(P')$, then their stabilizers $H_P$ and $H_{P'}$ are
conjugated.\\  
In this special case the Riemann-Hurwitz formula can be written as follows:  
$$2g(X)-2=|G|\Bigg(2g(X/G)-2+\sum _{i=1}^{s}\bigg(1-\frac{1}{v_{i}}\bigg)\Bigg)$$\\
[3 mm]
Applying this formula one can prove the following theorem, due to Kulkarni
 \cite{b5}.
We recall that a subgroup $G$ of $Aut(X)$ is said to be {\it large} if:
 $$|G|>4(g(X)-1)$$
\begin{theorem}[Kulkarni]\label{kulkarni}
Let $X$ be a Riemann surface of genus $g(X)\geq 2$, 
let $G$ be a large group of automorphisms of $X$ and consider the
quotient surface $Y:=X/G$.\\
Then $g(Y)=0$ and the projection map $\pi :X\to Y$ has 
either $3$ or $4$ critical values.\\
At most the following possibilities
  for the branching indices can occur.
\setlength{\tabbingsep}{0pt}
\begin{tabbing}
A) Four\=\ critical values\+\\[2mm]
a) \' One infinite family: \hspace{1cm}\= $\{ 2,2,2,n\}$\ \ \ $3$ \= $\leq n$\\
b) \' Three  exceptional cases:\> $\{ 2,2,3,a\}$  \>$3$ \' $\leq a\leq5$\-\\[3mm]
B) Three critical values\+\\[2mm]
a) \' One infinite family depending on $2$ parameters:\ \ $\{
2,m,n\}$\ \ \ $3\leq m\leq n$\\ 
 and if $m=3$ then $n\geq 7$, if $m=4$ then $n\geq 5$\\[1mm]
b) \' Five \=  infinite families:\+\\
i) \' $\{ 3,a,n\}$\ \ \ \ \ \ $3$ \=$\leq a$ \=$\leq6$,\ \ $a\leq n$; if $a=3$ then $n\geq 4$\\
ii) \' $\{ 4,4,n\}$  \>$4$ \'$\leq n$\-\\[1mm]
c) \' Exceptional cases:\+\\
i) \' $\{ 3,7,a\}$ \>$7$ \'$\leq a$\>$\leq 41$\\
ii) \' $\{ 3,8,a\}$ \>$8$ \'$\leq a$\>$\leq 23$\\
iii) \' $\{ 3,9,a\}$ \>$9$ \'$\leq a$\>$\leq 17$\\
iv) \' $\{ 3,10,a\}$ \>$10$ \'$\leq a$\>$\leq 14$\\
v) \' $\{ 3,11,a\}$ \>$11$ \'$\leq a$\>$\leq 13$\\
vi) \' $\{ 4,5,a\}$ \>$5$ \'$\leq a$\>$\leq 19$\\
vii) \'$\{ 4,6,a\}$ \>$6$ \'$\leq a$\>$\leq 11$\\
viii) \'$\{ 4,7,a\}$ \>$7$ \'$\leq a$\>$\leq 9$\\
ix) \' $\{ 5,5,a\}$ \>$5$ \'$\leq a$\>$\leq 9$\\
x) \' $\{ 5,6,a\}$ \>$6$ \'$\leq a$\>$\leq 7$
\end{tabbing}
\end{theorem}
\noindent In the table above we write the ramification 
indices of the critical values as a list:
    therefore, for example, $\{2,2,2,n\}$ means ``$4$ critical values,
      $3$ of index $2$ and the fourth of index $n$'', and $\{2,m,n\}$
        means ``$3$ critical values, one of index $2$, one of index $m$ and
	  the third of index $n$''\\[3mm]
So we have seven infinite families (one depending on two parameters) 
and $102$ exceptional cases.
\section{The abelian case}
\noindent In this section we assume $G$ abelian.\\
Let then $X$ and $Y$ be Riemann surfaces and let 
$\pi:X\to Y$ be a finite $G$-abelian cover, i.e. a Galois cover with finite 
abelian Galois group $G$. For a point $Q\in Y$ fix a point $P\in\pi^{-1}(Q)$,
then the stabilizer subgroup of $Q$
$$H_Q=\{\sigma\in G|\sigma P=P\}$$
is well-defined, since two subgroups of an abelian group are
conjugated if and only if they coincide.\\
Assuming $G$ abelian not all cases of Kulkarni's list do occur, but we have the following:
\begin{theorem}\label{listaclelia}
Let $X$ be a Riemann surface and $G$ a large abelian group of 
automorphisms of $X$.\\   
Then at most the following possibilities for the branching 
indices occur:
\begin{tabbing}
A) Four\=\ criti\=cal values\\[2mm]
\>\>$G\cong\mathbb{Z}_6$\hspace{2cm}\=$\{ 2,2,3,3\}$\\[3mm]
B) Three critical values\+\\[2mm]
a) \'  Three infinite families\+\\
i) \'  $G\cong \mathbb{Z}_{4g+2}$\>$\{ 2,2g+1,4g+2\}$\\
ii) \' $G\cong \mathbb{Z}_{4g}$\>$\{ 2,4g,4g\}$\\
iii) \'  $G\cong \mathbb{Z}_2\times \mathbb{Z}_{2g+2}$\>$\{ 2,2g+2,2g+2\}$\-\\[1mm]
b) \'  Ten exceptional cases.\+\\
i) \' $G\cong \mathbb{Z}_{12}$\>$\{3,4,12\}$\\
ii) \'  $G\cong \mathbb{Z}_{15}$\>$\{3,5,15\}$\\
iii) \'  $G\cong \mathbb{Z}_{6}$\>$\{3,6,6\}$\\
iv) \'  $G\cong \mathbb{Z}_{21}$\>$\{3,7,21\}$\\
v) \'  $G\cong\mathbb{Z}_{9}$\>$\{3,9,9\}$\\
vi) \' $G\cong \mathbb{Z}_{5}$\>$\{5,5,5\}$\\
vii) \' $G\cong \mathbb{Z}_3\times \mathbb{Z}_9$\>$\{3,9,9\}$\\
viii) \'  $G\cong \mathbb{Z}_3\times \mathbb{Z}_6$\>$\{3,6,6\}$\\
ix) \' $G\cong \mathbb{Z}_5\times \mathbb{Z}_5$\>$\{5,5,5\}$\\
x) \'  $G\cong \mathbb{Z}_4\times \mathbb{Z}_4$\>$\{4,4,4\} $
\end{tabbing}
\end{theorem}
\noindent We note that we are left with only three infinite families and $11$
exceptional cases.\\
In order to prove the previous theorem we
use the following result due to Nakajima \cite{b4}.
\begin{theorem}\label{nak}
Let $X\stackrel{\pi}{\to}Y$ be a finite $G$-abelian cover and let
$T_1, \ldots ,T_s\in Y$ be the critical values of $\pi$. For all
$j_0\in\{1,\ldots,s\}$  we denote by $G_{j_0}$ the subgroup of $G$
generated by all 
the stabilizers subgroup $H_{T_1},\ldots,H_{T_s}$ except $H_{T_{j_0}}$.\\
Then, if $g(Y)=0$ we have $G=G_{j_0}$.\\
Consequently, if $e_{T_{j}}:=|H_{T_j}|$ and $n_{j_0}:=\prod_{j\neq
  j_0}e_{T_j}$ then $n_{j_0}$ is a multiple of $|G|$.
\end{theorem}
\noindent From Nakajima's theorem it follows easily that the
\begin{corollary}\label{cor}
If $G$ is cyclic, then the least common multiple of any $s-1$ branching 
indices equals the order of $G$.
\end{corollary} 
\noindent Applying these results we can exclude many cases from the original
Kulkarni's list and prove the theorem \ref{listaclelia}.
\begin{prf}{listaclelia}
We start by considering the cases in Kulkarni's list with $4$ 
critical values.\\
$a)$ $\{ 2,2,2,n\},$\ \ \ $3\leq n$.\\
By theorem \ref{nak} $G$ must be generated by elements of order $2$,
contradicting the fact that it contains a cyclic subgroup of order $n \geq
3$. Therefore this case does not occur.\\
$b)$ $\{ 2,2,3,a\},$\ \ \ $3\leq a\leq 5.$\\
Let us start with the case $a=3$. In this case clearly $2\big| |G|$
and $3\big| |G|$ $\Rightarrow 6\big| |G|$. Moreover by theorem
\ref{nak} $|G|\big| 12$ and $|G|\big| 18$ $\Rightarrow |G| \big|
6$. Then $G\cong \mathbb{Z}_{6}$.\\
The same argument for $a=4$ gives  $12 \big| |G|$ and $|G| \big| 16$,
and for $a=5$ gives  $30 \big| |G|$ and $|G| \big| 12$, therefore both
cases do not occur.\\
We are left with the cases with $3$ critical values.\\
$a)$ $\{ 2,m,n\},$\ \ \ $3\leq m\leq n$.\\
In this case $2\big| |G|$, $m\big| |G|$ and $n \big| |G|$. By theorem
\ref{nak}, $|G| \big| 2m$, $|G| \big| 2n$ and $|G| \big| mn$. It follows
that either $|G|=n=m$ or  $|G|=n=2m$ or $|G|=2n=2m$.\\
i) $|G|=n=2m$. In this case $m \neq 3$ (by the condition $m=3
\Rightarrow n \geq 7$). Moreover (by $|G|=n$) $G$ is cyclic and therefore
we can use corollary \ref{cor} and conclude that $m$ is odd. By
Riemann-Hurwitz formula $m=2g_X+1$ and we have group $G\cong
\mathbb{Z}_{4g_{X}+2}$ and ramification indices
$\{2,2g_{X}+1,4g_{X}+2\}$.\\ 
ii) $|G|=2n=2m$. In this case $|G| \big| nm \Rightarrow m$ even.   
By Riemann-Hurwitz formula $m=2(g_X+1)$, $m \neq 4$ (by the condition
$m=4 \Rightarrow n \geq 5$). Moreover, by corollary \ref{cor}, $G$
cannot be cyclic and therefore $G\cong
\mathbb{Z}_{2} \times \mathbb{Z}_{2g_{X}+2}$ and the ramification indices
are $\{2,2g_{X}+2,2g_{X}+2\}$.\\ 
iii) $|G|=n=m$. Arguing as in the previous cases $G$ is cyclic and
$4 \neq m=4g_X$. Therefore $G\cong \mathbb{Z}_{4g_X}$ and the
ramification indices are $\{2,4g_{X},4g_{X}\}$.\\ 
$b)$ $\{ 3,a,n\},$\ \ \ $3\leq a\leq 6$ and $\{ 4,4,n\}$.\\
If $a=3$, Kulkarni's list imposes $n \geq 4$. By $n \big| |G| \big|
3a=9$, we conclude $n=|G|=9$ and therefore the group is cyclic and we
can apply corollary \ref{cor} and get a contradiction.\\
If $a=4$ or $a=5$ we have $3a \big| |G| \big| 3a$ and therefore clearly
$G\cong \mathbb{Z}_{3a}$ and by theorem \ref{nak} we have ramification
indices $\{3,a,3a\}$.\\
If $a=6$ we have $6 \big| |G| \big| 18$. If $|G|=6$, $G \cong
\mathbb{Z}_{6}$ and the ramification indices are
$\{3,6,6\}$. Otherwise $|G|=18$ and (since by corollary \ref{cor}
$G$ is not cyclic) the only possibility is $G\cong \mathbb{Z}_{3}
\times \mathbb{Z}_{6}$ and ramification indices $\{3,6,6\}$.\\
For the case $\{4,4,n\}$, $n \geq 4$ we first note that by theorem
\ref{nak} $G$ is generated by two cyclic subgroups of order $4$ and
therefore no element of $G$ has higher order. In particular $n=4$.  
We have $4 \big| |G| \big| 16$. If $|G|=4$ by Riemann-Hurwitz formula
$$ 2g_{X}-2=4\bigg(1-\frac {3}{4}\bigg)=1,$$
a contradiction.\\
Then $|G|= 8$ or $16$ and $G$  contains three cyclic
subgroups of order $4$ such that each pair of them generates $G$: the
only possibility is $G \cong \mathbb{Z}_{4} \times \mathbb{Z}_{4}$.\\ 
Note that we have already found all the cases in our
list except the $5$ ``exceptional'' that we have labeled as iv), v),
vi), vii) and ix). 
It remains to check the last 99 exceptional cases (case ``c)'') of
Kulkarni's list and show that (for $G$ abelian) only these $5$ are
possible: this can be proved with the same arguments we have already used
in this proof, and  we leave the details to the reader.
\end{prf}
\section{Constructing curves with large abelian group of automorphism:
the main theorem}
\noindent In the previous section we have restricted Kulkarni's list in the abelian
case to  three infinite families and $11$ exceptional cases; in this
section we show that all these cases do occur, constructing
the curves. 
\begin{theorem}\label{main}
The Riemann surfaces with a large abelian group of automorphisms and
$4$ critical values are desingularisation of the $1-$parameter family
of curves of genus $2$ in the weighted projective space
$\mathbb{P}(2,1,1)$ of equation 
$$y^6=x^3z^3(x-z)^2(x-\lambda z)^4$$
where $y$ is the coordinate of weight $2$ and $\lambda$ varies in
$\mathbb{C}\setminus\{0,1,-1\}$.\\ 
The group of these curves is isomorphic to
$\mathbb{Z}_6$, and the ramification indices are $\{2,2,3,3\}$.\\
The other Riemann surfaces with a large abelian group of automorphisms
are listed in the following table:
\begin{center}
\setlength{\extrarowheight}{2pt}
\begin{tabular}{|c|c|c|c|r@{\hspace{3pt}}c@{\hspace{3pt}}l|} \hline
$1$& $\Z_{4g+2}$ & $2,2g+1,4g+2$ & $g\geq 2$ & $y^{4g+2}$ & $=$ & $x^{2g+1}(x-z)^2 z^{2g-1}$
\\[2pt] \hline
$2$ & $\Z_{4g}$ & $2,4g,4g$ & $g\geq 2$ & $y^{4g}$ & $=$ & $x^{2g}(x-z)^2 z^{2g-1}$\\[2pt]
\hline
$3$ & $\Z_{12}$ & $3,4,12$ & $3$ & $y^{12}$ & $=$ & $x^{4}(x-z)^3 z^{5}$\\[2pt] \hline
$4$ & $\Z_{15}$ & $3,5,15$ & $4$ & $y^{15}$ & $=$ & $x^{5}(x-z)^3 z^{7}$\\[2pt] \hline
$5$ & $\Z_{6}$ & $3,6,6$ & $2$ & $y^{6}$ & $=$ & $x^{4}(x-z) z$\\[2pt] \hline
$6$ & $\Z_{21}$ & $3,7,21$ & $6$ & $y^{21}$ & $=$ & $x^{7}(x-z)^3 z^{11}$\\[2pt] \hline
$7$ & $\Z_{9}$ & $3,9,9$ & $3$ & $y^{9}$ & $=$ & $x^{3}(x-z)^5 z$\\[2pt] \hline
$8$ & $\Z_{5}$ & $5,5,5$ & $2$ & $y^{5}$ & $=$ & $x(x-z) z^{3}$\\[2pt] \hline
$9$ & $\Z_2\times\Z_{2g+2}$ & $2,2g+2,2g+2$ & $g\geq 2$ & 
{\begin{tabular}{r@{}}
$y^2$\\
$w^{2g+2}$
\end{tabular}} & 
{\begin{tabular}{@{}c@{}}
$=$\\
$=$
\end{tabular}} &
{\begin{tabular}{@{}l@{}}
$xz$\\
$(x-z)z^{2g+1}$
\end{tabular}}\\[3pt] \hline
$10$ & $\Z_3\times\Z_{9}$ & $3,9,9$ & $7$ &
{\begin{tabular}{r@{}}
$y^3$\\
$w^{9}$
\end{tabular}} &
{\begin{tabular}{@{}c@{}}
$=$\\
$=$
\end{tabular}} &
{\begin{tabular}{@{}l@{}}
$xz^2$\\
$(x-z)z^{8}$
\end{tabular}}\\[3pt] \hline
$11$ & $\Z_3\times\Z_{6}$ & $3,6,6$ & $4$ &
{\begin{tabular}{r@{}}
$y^3$\\
$w^{6}$
\end{tabular}} &
{\begin{tabular}{@{}c@{}}
$=$\\
$=$
\end{tabular}} & 
{\begin{tabular}{@{}l@{}}
$xz^2$\\
$(x-z)z^{5}$
\end{tabular}}\\[3pt] \hline
$12$ & $\Z_5\times\Z_{5}$ & $5,5,5$ & $6$ & 
{\begin{tabular}{r@{}}
$y^5$\\
$w^{5}$
\end{tabular}} &
{\begin{tabular}{@{}c@{}}
$=$\\
$=$
\end{tabular}} &
{\begin{tabular}{@{}l@{}}
$xz^4$\\
$(x-z)z^{4}$
\end{tabular}}\\[3pt] \hline
$13$ & $\Z_4\times\Z_{4}$ & $4,4,4$ & $3$ & 
{\begin{tabular}{r@{}}
$y^4$\\
$w^{4}$
\end{tabular}} & 
{\begin{tabular}{@{}c@{}}
$=$\\
$=$
\end{tabular}} & 
{\begin{tabular}{@{}l}
$xz^3$\\
$(x-z)z^{3}$
\end{tabular}}\\[3pt] \hline
\end{tabular}
\end{center}
In the second column we wrote the group $G$, in the third the
three ramification indices, in the fourth the genus of the curve, in
the last the equations of a birational model of the curve in $\mathbb{P}^2$ or
in $\mathbb{P}^3$: in all cases the coordinates are chosen so that
the map $\pi$ is induced by the projection given by the coordinates $(x,z)$.
\end{theorem}
\noindent Note that there are two cases in the above list  such that the given
group admits a subgroup that is still large: the families ``1'' and
``9'', in the  case $g=2$, admit a subgroup of index $2$ that is still
large. Looking at the list one easily concludes that these are
respectively the cases ``8'' and ``5''. So, if one is only interested
in the list of Riemann surfaces having a large group of automorphisms,
one can drop the cases $5$ and $8$.\\
To prove the theorem we use a result of Pardini \cite{b3} that
classifies abelian covers; we will briefly state her (much more
general) result only in the special case of Riemann surfaces. \\ 
Let $X$, $Y$ be Riemann surfaces and let $\pi:X\to Y$ be a finite
abelian $G$-cover.\\ 
The direct image sheaf $\pi_{*}(\mathcal{O}_X)$ splits as a direct sum of
line bundles as follows: 
\[\pi_{*}({\mathcal O}_{X})=\bigoplus  _{\chi \in G^{*}}L_{\chi }^{-1}=
{\mathcal O}_{X}\oplus \Bigg(\bigoplus_{\chi \in G^*\setminus \{ 1\}}
L_{\chi}^{-1}\Bigg)\] 
where $G^*=Hom(G,\mathbb{C}^*)$ denotes the group of characters of $G$
and $G$ acts on $L_{\chi }^{-1}$ via the character $\chi$.\\
Let $Q_1,\ldots,Q_s$ be the critical values of $\pi$, and let us
associate to each critical value $Q_i$ its ramification index $e_i$,
its stabilizer subgroup $H_i$ and a generator $\psi_i$ of $H_i^*$.\\ 
Let $\chi_1,\ldots ,\chi_n\in G^*$ be such that $G^*$ is the direct
sum of the cyclic subgroups generated by the $\chi_j$'s, let $d_j$ be
the order of $\chi_j$ and write (for the sake of simplicity) $L_j$ for
$L_{\chi_j}$.\\  
To each pair $Q_i, \chi_j$ we associate $r_{i,j}$ such that 
$$ (\chi_j)_{|H_i}=\psi_i^{r_{i,j}} \ \ \ \ 0\leq r_{i,j}<e_i.$$
We have the following:
\begin{theorem}\label{parda}
Let $\pi \colon X\to Y$ be a finite abelian $G$-cover of Riemann surfaces.\\
Then we have the following set of linear equivalences:
$$d_jL_j\equiv \sum _{i=1}^s\frac {d_jr_{i,j}}{e_i}Q_i \ \ \ \ j=1,\ldots ,n$$
Conversely, for assigned $G$, $Q_i$, $H_i$, $\psi_i$, $i=1,\ldots,s$ to each 
set of line bundles $L_j$, $j=1,\ldots ,n,$ satisfying the previous linear
equivalences there corresponds a unique, up to isomorphism, 
abelian $G$-cover of $Y$ branched on $Q_1,\ldots,Q_s$ and such that 
$H_i$ is the stabilizer subgroup of $Q_i$ and $\psi_i$ the corresponding
character.
\end{theorem}
\begin{prf}{main}
Let us first write the linear equivalences of theorem \ref{parda} in
the case of cyclic covers of $\mathbb{P}^1$. In this case there is
only one linear equivalence.\\
Let $Q_1,\ldots ,Q_s\in \mathbb{P}^1$ be the critical values,
$H_1,\ldots ,H_s$ and $e_1,\ldots,e_s$ the respective stabilizer
subgroups and ramification indices.\\
For every generator $\chi$ of $G^*$ there is a generator $g$ of $G$
such that 
$$\chi(g)=e^{\frac{2\pi i}{n}},$$
where $n$ is the order of $G$.
$H_i$ is generated by $g^{\frac{n}{e_i}}$.\\
For every generator $\psi_i$ of $H_i^*$ there is an $\alpha$ with
$(\alpha,e_i)=1$ (so that $g^{\alpha \frac{n}{e_i}}$ generates $H_i$)
such that $$\psi_i\Big(g^{\alpha \frac{n}{e_i}}\Big)=e^{\frac{2\pi
    i}{e_i}}.$$ 
It follows that
$$\chi_{|H_i}\Big(g^{\alpha \frac{n}{e_i}}\Big)=e^{\frac{2\pi i}{n}\cdot \alpha \cdot\frac{n}{e_i}}=
e^{\frac{2\pi i}{e_i}\cdot \alpha} $$
We have then to solve the linear equivalence
$$nL\equiv \sum _{i=1}^s\frac{n}{e_i}\alpha_i Q_i$$
for $(\alpha_i,e_i)=1$,\ \ \ $0\leq \alpha_i<e_i$.\\
Since two divisors on $\mathbb{P}^1$ are linearly equivalent if and only 
if they have the same degree, we only have to find solutions
$\alpha_i$ of the congruence 
$$ \sum _{i=1}^s\frac{n}{e_i}\alpha_i\equiv 0\ \ (\mbox{mod}\ n)$$
with  $(\alpha_i,e_i)=1$, $0\leq \alpha_i<e_i$.\\
We can now use theorem \ref{parda} to construct all the pairs $X,G$
in the list of theorem \ref{listaclelia} with $G$ cyclic.\\
Consider the case:\\
 $G\cong \mathbb{Z}_{4g+2}$\ \ \ $\{2,2g+1,4g+2\}$\\
Let
$Q_{1},\ Q_{2},\ Q_{3}$ be the critical values and $H_{1},\ H_{2},\ H_{3}$
the respective stabilizers.\\
By theorem \ref{nak} $G$ is generated by each pair of stabilizers; so
$G\cong  H_{1}\times H_{2}$ with $H_{1}\cong \mathbb{Z}_{2}$,
$H_{2}\cong \mathbb{Z}_{2g+1}$.\\
We have also an isomorphism
$$\begin{array}{ccc}
G^{*}&\stackrel{\cong}{\longrightarrow}& H_{1}^{*}\times H_{2}^{*}\\
\chi& \mapsto& (\chi_{|H_{1}},\chi_{|H_{2}})
\end{array}
$$
In particular for each pair $(\psi_{1},\psi_{2})$ such that
$\langle \psi_{1}\rangle =H_{1}^{*}$ and $\langle \psi_{2}\rangle =H_{2}^{*}$
there is a unique generator $\chi$ of $G^{*}$
such that $$\chi_{|H_{1}}=\psi_{1}\ \ \ \chi_{|H_{2}}=\psi_{2}$$
The linear equivalence to be solved is then:
$$(4g+2)L\equiv(2g+1)Q_{1}+2Q_{2}+\alpha Q_{3}$$
with $0<\alpha<4g+2,\ \ (\alpha,4g+2)=1$.\\
The only solution is $\mbox{deg}L=1$, $\alpha=2g-1$.
This shows the existence of this case.\\ 
In fact this method gives these curves explicitely, as desingularisation of 
the curves in $\mathbb{P}^2$ of equation 
$$y^{4g+2}=x^{2g+1}(x-z)^2 z^{2g-1}.$$
The large group of automorphisms of these curves is isomorphic to
$\mathbb{Z}_{4g+2}$, generated by the automorphism induced on the curve by the
following automorphism of $\mathbb{P}^2$ of order $4g+2$:
$$(x,y,z)\mapsto (x,e^{\frac{2\pi i}{4g+2}}y,z)$$
The projection to the quotient $\pi$ is induced by the projection of
center $(0,1,0)$, the coordinates $(x,z)$ of the $3$ critical values
of $\pi$ are $Q_1=(1,0)$, $Q_2=(1,1)$ and $Q_3=(0,1)$, 
of respective ramification index $2,2g+1$ and $4g+2$.\\
The other cases with cyclic group and $3$ critical values are
similar, and we leave the details to the reader.\\ 
The one with $4$ critical values is slightly different:\\
$G\cong \mathbb{Z}_6$\ \ \ $\{ 2,2,3,3\}$\\
In this case $G\cong H_1\times H_3$ and analougsly 
$G^{*}\cong H_1^{*}\times H_3^{*}$.\\ 
The linear equivalence to be solved is
$$6L\equiv 3Q_1+ 3 Q_2+2Q_3 +2\beta Q_4$$
that has the only solution $\mbox{deg}L=2$, $\beta =2$.\\
The difference is in the explicit construction of the curves, since
$\deg L \neq 1$. We have a one-parameter family of
curves of genus $2$ that are desingularisation of the curves in the
weighted projective space $\mathbb{P}(2,1,1)$ of equations
$$y^6=x^3z^3(x-z)^2(x-\lambda z)^4$$
where $y$ is the coordinate of weight $2$ and $\lambda$ varies in
$\mathbb{C}$.\\ 
We are left with the cases where $G$ is not cyclic. Consider for example\\ 
$G\cong \mathbb{Z}_2 \times \mathbb{Z}_{2g+2}$\ \ \ $\{2,2g+2,2g+2\}$\\
Since each pair of stabilizers generates $G$, $G \cong
H_1\times H_2$ with $H_1\cong \mathbb{Z}_2$ e $H_2\cong
\mathbb{Z}_{2g+1}$. \\
Analougsly $G^* \cong H_1^* \times H_2^*$: for each pair
$(\psi_1,\psi_2)$ with $\langle \psi _1 \rangle =H_1^{*}$ and
$\langle \psi _2 \rangle =H_2^{*}$, there are two uniquely determined
characters $\chi_1$ and $\chi_2$ such that $\langle \chi_1,\chi_2
\rangle =G^{*}$ and 
\begin{center}
\begin{tabular}{ll}
$\chi _{1|H_1}=\psi _1$\hspace{2cm}\ &$\chi_{2|H_1}=1$\\
$\chi _{1|H_2}=1$&$\chi_{2|H_2}=\psi _2$
\end{tabular}
\end{center}
By theorem \ref{parda} we have to solve two linear equivalences:
$$\left\{
\begin{array}{ll}
2L_1\equiv Q_1+ \frac {\alpha}{g+1} Q_3& \text{ with }\ \  
0<\alpha<2g+2,\ \ \text{ therefore }\alpha\in\{0,g+1\}\\
(2g+2)L_2\equiv Q_2+ \beta Q_3 &\text{ with }\ \  
0<\beta <2g+2
\end{array} \right .
$$
The only solution is $\deg L_1=\deg L_2=1$, $\alpha =g+1,\beta=
2g+1$. This shows the existence of this case.\\
These curves are the desingularisation of the curve of equations
in $\mathbb{P}^3$
$$\left\{ 
\begin{array}{l}
y^2=xz\\
w^{2g+2}=(x-z)z^{2g+1}
\end{array}
\right .
$$
The projection to the quotient is the one induced by the projection 
$\mathbb{P}^3 \dashrightarrow \mathbb{P}^1$ of center $\{x=z=0\}$. $G$
is generated by 
$$\left\{ 
\begin{array}{lll}
(x,y,z,w) & \mapsto& (x,-y,z,w) \\
(x,y,z,w) & \mapsto& (x,y,z,e^{\frac{2\pi i}{2g+2}}w) 
\end{array}
\right .
$$
The same method works in all cases; we leave the details to the reader.
\end{prf}
We have the following
\begin{corollary}\label{ultimo}  
Let $X$ be a Riemann surface and suppose that on $X$ acts an abelian large 
group of automorphisms. If $g(X)\geq8$ then $X$ is hyperelliptic.\\
More precisely all the curves in the theorem \ref{main} are
hyperelliptic except the cases 
$3$, $4$, $6$, $7$, $10$, $11$, $12$, $13$.
\end{corollary}
\begin{prf}{ultimo}
The family of Riemann surfaces with $4$ critical values and the cases
$5$ and $8$ have genus $2$ and therefore are hyperelliptic.\\
In the cases $1$, $2$ and $9$ it is easy to find an involution in $G$
with $2g+2$ fixed points.\\
It remains to show that all the remaining cases are not hyperelliptic.\\
In fact, if the hyperelliptic involution $\tau$ exists, it belongs to the
center of the automorphisms group. It follows that if $\tau\notin G $,
there exists  $G'$ abelian subgroup of $Aut(G)$ with $|G'|=2|G| \geq
2(4g-3)$, contradicting (since $g\geq3$) the upper bound $4g+4$ for
the cardinality of an abelian group of automorphism. Therefore, if
$\tau$ exists, $\tau \in G$.\\ 
It is now easy to conclude: the groups in the cases $4$, $6$, $7$, 
$10$ e $12$ do not contain any involution, and one can check that the
involutions in cases $2$, $11$ and $13$ have not $2g+2$ fixed points.
\end{prf}

\begin{thebibliography}{10}
\bibitem{hur} Hurwitz, A.,
{\it \"Uber algebraische Gebilde mit eindeutige Transformationen in
  sich}.
Math.~Ann.~41 (1893), 403--442.

\bibitem{klein}
Klein, F.,
{\it \"Uber die Transformationen siebenter Ordnung der elliptischen
  Funktionen}.
Math.~Ann.~14 (1879), 428--471.

\bibitem{b5} 
Kulkarni, R.~S., 
{\it Riemann surfaces admitting large automorphism groups}.  
Extremal Riemann surfaces (San Francisco, CA, 1995),  63--79,
Contemp.~Math., 201, Amer. Math. Soc., Providence, RI, 1997. 

\bibitem{b2} Miranda, R.,
{\it Algebraic Curves and Riemann Surfaces}.
Graduate Studies in Mathematics, 5, Amer. Math. Soc., 1995.

\bibitem{b4} 
Nakajima, S.,
{\it On abelian automorphism groups of algebraic curves}.  
J.~London Math.~Soc.~(2)  36  (1987),  no. 1, 23--32.

\bibitem{b3} 
Pardini, R., 
{\it Abelian covers of algebraic varieties}. 
J.~Reine Angew.~Math.\  417  (1991), 191--213.

\bibitem{schwartz}
Schwartz, H.\ A.,
{\it \"Uber diejenigen algebraischen Gleichungen zwischen zwei
  ver\"anderlichen Gr\"ossen, welche eine Schaar rationaler, eindeutig
  umkehrbarer Transformationen in sich selbst zulassen}. 
J.~Reine Angew.~Math. 87 (1890), 139--145.

\bibitem{wiman}
 Wiman, A.,
{\it \"Uber die hyperelliptischen Kurven und diejenigen vom Geschlecht
  $p=3$, welche eindeutige Transformationen in sich zulassen}.
 Bihang Till.~Kungl.~Svenska Vetenskaps-Akademiens Handlingar 21,
 no. 1, (1895), 1--23.


\end{thebibliography}

\end{document}